# Sustainable Closed-loop supply chain under uncertainty

[1]Komeil Baghizadeh

[1]Department of industrial engineering, Khrazmi University, Tehran, 4681973333, Iran.

Corresponding author: Komeil Baghizadeh; std_komeilbaghizadeh@khu.ac.ir ;
Komeil.baghizadeh@gmail.com .

## Introduction

Currently, sustainable economic models are important tools for companies in supplying natural resources sustainably due to rising cost pressures with economic sustainability, environmental and political constraints in the supply chain, and environmental challenges,.

Recently, a large body of research has been devoted to supply chains to decline the material and source burden of fabrication, products, and waste via resource efficacy (Daaboul et al. 2016; Genovese et al. 2015), and green and low-carbon supply chain enterprises (Pan et al. 2015; Park, Sarkis, and Wu 2010; Zhu, Geng, and Lai 2010); however, environmental sustainability and economic performance tend to be seen as a tradeoff (Colicchia et al. 2016). Thus, a continuous challenging task is to retain competitiveness and concurrently create social and environmental value via supply chain remodel. The framework and construct of discourse can be differently established concerning value creation feasibilities from the remodel of supply chains and seepages of product and materials structurally, which arise from geographical dispersal and complicated multifaceted rules of materials and product intricacy (Klassen 2009). A strategy is closed loop supply chain (CLSC) plan described by Guide and Van Wassenhove (2009, 10) as 'designing, controlling, and operating in a system for maximizing value creation during the whole lifespan of a product to dynamically recover value from various kinds and bulks of revenues during de years, and comprises managing product return, leasing, and refabricating (Blackburn et al. 2004; Hu et al. 2014; Klassen 2009).

Recently, academically published research has progressively focused on the term circular economy (CE) in addition to the CLSC (e.g. Bocken et al. 2016; Ghisellini, Cialani, and Ulgiati 2016). A formally defined CE based on the present study is the one 'designed restoratively and renewably , aiming at keeping products, constituents, and materials at their uppermost efficacy and value at any moment to distinguish between technical and biological cycles'.1 The design of such an economy aimed at preserving and enhancing natural wealth, optimizing resource profits, and minimizing system risks through management of limited stocks and recoverable flows (Webster 2013, 2015)

An increasing agreement exists on switching from our current industrial linear model to a CE model as the only approach ahead with sustainable fabrication and development (Geissdoerfer et al., 2017; Su et al., 2013)

A novel economic model is progressively emerging because the linear economy model's frequently fail to fulfil the dynamic sustainability issues worldwide. The concept of CE has received recognition from governments, companies, and universities. CE is essentially based on the notion that products, production systems, and supply chains need to be transformed for the establishment of practical associations between ecological systems and economic growth, which also pushes the boundaries of environmental sustainability. It focuses on the construction of self-sustaining production systems where materials are reutilized repeatedly (Genovese et al., 2015). It is necessary to incorporate these

CE principles into the designing and managing policies for supply chains in energy systems to minimize material flows and reduce accidental undesirable fallouts of production procedures (Srivastava, 2007).

As the linear economy transitioned to CE, businesses should take sustainability and closed-loop cycle into consideration. The emphasis of CE is on minimum resource consumption and environmental protection, therefore businesses seek to implement green supply chain. Wastes are particularly produced by traditional supply chain that lead to ecological complications and the society and environment are not taken into account.

Discarding consumer goods is prevented through the closed-circuit procedures resulting from environmental and social attentions. Inside this closed-circuit, the CE puts forward a system-centered strategy regarding the recycle and reuse of materials and to diminish the volume of materials, thereby making sustainability more promising (Sauv?e et al., 2016).

The Ellen MacArthur Foundation [1] as the chief paragon beside other entities has realized that the present globe simultaneously needs persistent economic development, environmental protection, and societal welfare. Most often, however, the CE notion that is pressed in politically related areas suffers the lack of a demarcated technological or industrial application. wherein this situation, the 6R approach (Reduce, Reuse, Recycle, Recover, Redesign, and Remanufacture) launched for sustainable engineering in recent years and novel breakthroughs should be required for creating and defining an innovative commercial model, resulting in sustained value creation for designing and fabricating products [2].

The CE is achievable via protracted designing, preservation, repairing, reprocess, remanufacture, renovating, and recovering [1]. If these novel accomplishments are extensively introduced into existing industrial systems, it will possibly be environmentally advantageous positively whereas it disrupts the working manner of current entities and supply chains, which necessitates redesigning existing supply chains.

The idea of CLSCs has emerged to respond to the requirement for reinterpreting supply chains in the CE framework. Investigational work has often focused on closed loops at the end-of-life (EOL) of consumptive products. Nonetheless, the loop is closed throughout product lifespans and for the rest of principal and ancillary product types. It is interesting that only scarce research work tries to discriminate various kinds of closed loops in supply chains. For instance, Wells and Seitz [2] proposed four kinds of closed loops—internal or within fabrication, post-business, post-consumer, and post-society—and discussed a variety of features and disadvantages for every loop. Such a restricted examination of closed loops calls for additional investigation on the issue, in particular because novel CE applications and commercial models become popular in companies. Although CLSCs are pivotal to CE business models, they have received research work principally as discrete investigational topics. Additionally, little studies are available on clarifying the association of CE and sustainable supply chain management (SSCM) publications and applications.

A notion helping supply chain within the framework of CE is industrial symbiosis with the aim of abstracting maximal profits of resources, products, energy consumption, and synergizing it to achieve further sustainability throughout the supply chain [3, 4]. In the framework of industrial symbiosis, a paper-fabricating institute was surveyed by focusing on 6Rs, and a discussion was made on the existing opportunities with modern technologies including Internet of Things (IoT) and Artificial Intelligence (AI) in Industry 4.0 era [5]. Industry 4.0 integrates the physical and digital world to unlock the digitalization advantages to promote the sustainability of the institute and in turn help to develop CE

[6, 7]. Industry 4.0 ensures detection of quality issues at the onset, which enhances the efficiency and declines waste expenditure considerably [9, 10]. The sustainable supply chain, and concurrently 4.0, can be described as the supply chain that utilizes Industry 4.0 instruments with the aim of closing the material and energy cycles, in addition to assisting the information flow and the activities, yielding more effective, smarter, more accurate, and faster operations.

in Germany, the manufacturing industry proposes and adopts Industry 4.0 as part of the High-Tech Strategy 2020 Action Plan [1] through implementation of an intelligent factory on the basis of linking and observation.

I4.0 technologies can be advantageous to managing operation in a variety of means, including reducing the process time of a product, decreasing manufacture cost, up-grading the value chain coordination, improving processing resilience, enhancing customer service, greater product customization, etc. (Fettermann et al., 2018). Four maturity levels recommended by Fettermann et al. (2018) include (i) technologies, (ii) control, (iii) optimization, and (iv) autonomy. The technologies comprise instantaneous information regarding observing and reporting at the interior level (e.g. operative units including inventory-dependent parameters, set-up time, etc.) and exterior level (e.g. demand, lead time, order dimension, transportation management, etc.). The system of operations is managed by application software at the interior and exterior levels of services (e.g. cloud enterprise resource planning (ERP) system outfitted with an inventory management module). The level optimization involves using algorithms capable of optimizing the system operations. Ultimately, a system is adaptable to the environment to promote performance at the autonomic level.

ICPT (Information, Communication, and Production Technology) can expectedly be applied as an instrument in the future generation of engineering setting including intelligent factories and encompasses the Internet of Things (IoT), big data analytics, cloud computing, the Cyber-Physical System (CPS), Artificial Intelligence (AI), Virtual and Augmented Reality (VR & AR), and 3D printing technology. These seem to be congregated in the industrial sector to create new values that differ with traditional ones. Besides, ICPT is denoted as a general technology with capabilities of visibility, flexibility, responsiveness, integrity, and automaticity for a novel manufacture prototype.

Additionally, the conjunction of ICPT has the tendency to be incorporated not only into the manufacture industry, but also in the logistics and wholesale industries. An intelligent factory utilizing ICPT is capable of managing information and possesses multipurpose constituents. In logistics, ICPT is known as smart logistics [5,6], and retail as smart retail [7,8]. Intelligent logistics and retail joined with ICPT can track products and stock list in the supply chain, obtain prompt information of operations, promote decision-making challenges, and interplay with clients for a custom-made service to guarantee the competitive ability of their companies.

The intelligent factory that uses ICPT is regarded as a cutting-edge factory to promote throughput and effectiveness in resource usage and allotment, scheduling, control, optimization, etc. The processing time a factory can be improved by the use of a realistic decision-making procedure [15]. Shellshear et al. [16] denoted different topics concerning the association of data and cooperation amongst the constituents in a smart factory. These can promote efficacy in using resources and the flexible nature of the manufacture system [17] for optimizing companion choice with cooperation [18], and designing a strategical program and system with data association and its application in a realistic manufacture setting [19].

## Literature

**Supply chain based on circular economy**

The CE necessitates a novel association with products and materials, which is more labor intensive and less resource demanding [4]. Lacy and Rutqvist [5] discovered the way businesses can gain profit within the CE by 'value creation from wastage'. They recognized 'Recovery and Recycling' as a kind of business model within the CE, where all things being commonly regarded as wastage is recovered for other usages. They recognized two variable parameters within this model: the retrieval of value from EOL products and the retrieval of wastage and byproducts from a firm's own production procedure and practices. The second parameter involves the notion of industrial symbiosis (IS) referring to the yield wastage and underused flows of one institute that is utilized as valued and prolific inputs for another institute [6, 7]. In an ordinary entity, the discovery of sellers and consumers of these materials is under the horizon of supply chain management (SCM). SCM is, therefore, a prerequisite for making an efficacious IS usage and association [8, 9].

They play a role in advancing the understanding of CLSCs through offering a series of closed-loop types within CE business models and developing a conceptuality of CLSCs being demonstrated by realistic examples. This procedure will eventually provide backing to businesses for designing the CLSCs and CE business models, which therefore facilitates the acceptance of these operations in existing companies.

Among academics, an attractive research topic has been the creation of innovative CE business models for many years. Accordingly, it is interesting to understand the development course of sound business models within the CE [12, 13]. To build business model approaches within the CE, three mechanisms have been recognized as the model basics [14]: (1) decelerating resource loops via expansion or intensification of the usage period of products, for example product lifespan expansion via repairing or remanufacture, (2) closing resource loops with improving recycle applications, for example post-consumer plastics recycle, and (3) thinning resource loops through declining the volume of resources per product.

To transform from a linear to a circular economy, organizations need to renovate their supply chain. Based on this viewpoint, the CE can effectively act in promoting from the traditional supply chain to the CLSC and extra burden on the company (Zhu et al., 2010).

Geissdoerfer (2018) describe Circular Supply Chain Management (CSCM) as the conformation and harmonization of the organizational practices, marketing, sales, R & D, production, logistics, IT, finance, and client service within and across business units and organizations for closing, slowing, strengthening, narrowing, and dematerializing material and energy loops to minimalize resource input into and wastage and emission seepage out of the system for improving its operational efficacy and generating competitive advantages.

As asserted by Webster (2015), CE aims to maintain the maximal level usefulness and worth of the products and materials through designing, maintaining, repairing, reusing, remanufacture, and recycle while Merli et al. (2018) included declining wastage in the description. Geissdoerfer et al. (2017) denoted that the CE refers to a renewing system such that resource input wastage, emission, and energy usage are minimalized by closed loops of material and energy. Zhu et al. (2011) announced that CE aims to minimize the use of material and emissions, thereby supporting eco-design, cleanser fabrication, and waste control.

Schenkel, Krikke et al. (2015) emphasized the absence of investigation on value creation by various loops practically, in spite of a progressive attention to CLSC. Schenkel, Krikke et al. (2015) recognized four kinds of value creation, viz. economic, environmental, information, and consumer-centered. Two leading publications regarding the positive business and economic cases for a CE (EMF 2012, 2013) encouraged extensive attention to circular value creation potential resulting in several scientific and strategic papers on ways of achieving this goal (Bakker et al. 2014; Ghisellini, Cialani, and Ulgiati 2016; Haas et al. 2015; Lacy and Rutqvist 2015; Lieder and Rashid 2016), which supports the opinions of Park, Sarkis, and Wu (2010) about the promising contribution of CE to creating value within SCM.

This viewpoint results in an integrated approach with the growth of supply chains by taking both forward and reverse supply chains into consideration at the same time, reputed as CLSC [17]. Reverse supply chain comprises activities the deal with value retrieval of EOL products either by the main product producer or a third party [15, 18, 19]. EOL products are gathered from clients, followed by performing suitable procedures, including repair, disassemble, remanufacture, recycle and disposal of the products with environmental sensitivity [19]. A terminological difference should be noted between the business model mechanisms and the supply chain publications. Indeed, CLSC publications denote repairing and remanufacture to be a mechanism for closing the loop of a certain product while the above CE business model publications consider repairing and remanufacture as part of 'slowing the loop'. This exemplifies the disintegration of the CE by reviewing across multiple investigational disciplines, as emphasized by De Angelis et al. [20], which necessitates further research for developing a communal understanding between disciplines. A group of researchers distinguish between open-loop and closed-loop SCM [15]. In the latter, a constituent will be reutilized or reprocessed for the same usage, while the materials or constituents enter a different usage in the former view [21, 22]. CLSCs serve to tak back products from clients and return them to the main producer for the retrieval of added value by the reuse of the product entirely or partly [19]. Open loop supply chain encompasses materials retrieved by parties different from the main manufacturers who are able to reuse these materials or products [21]. The line between closed-loop and open-loop methods is mostly very thin; besides, the chief aim is to recover added value and avoid wastage, which gains support from converse logistics actions [23].

CE assists managers in focusing on the economic, social, and environmental achievements in a supply chain framework (Hamprecht et al., 2005; Ross et al., 2012). CE-based sustainability enterprises help in reducing the food waste and effect of contamination and improving the general performance via a variety of R's (recycle, reuse, reduce, etc.) (Yong, 2007; Geng et al., 2013; Govindan and Hasanagic, 2018). CE-based sustainability facets also help provide food with safer and better quality to the clients (Beske et al., 2014).

Pan et al. (2015) performed an analysis of several wastes to energize technologies such as combustion, gasification, and anaerobic digestion for providing portfolio choices of technologies for various kinds of wastes to energize supply chains for generating a CE system. Likewise, Nasir et al. (2016) utilized a case study from the construction industry for demonstrating and comparing the environmental profits that are achievable via adopting CE principles compared with the conventional linear production systems. Ahn et al. (2015) designed a deterministic mathematical programming model for strategic planning design of a biomass-to-biodiesel supply chain network from feedstock fields to end users that can concurrently satisfy resource constraints, demand constraints, and technology during a protracted planning prospect. Genovese et al. (2015) made a comparison between the performances of conventional and circular production systems throughout a series of indicators by two case studies from chemical and waste food (waste cooking oil to biodiesel) supply chains. The authors drew a conclusion that the incorporation of CE principles into sustainable supply chain management practices offers distinct environmental benefits.

Calderon et al. (2017) introduced a general optimization platform according to a multi-period mixed integer linear programming model to scrutinize the strategic designing of waste to synthetic natural gas supply chains.

Yılmaz Balaman (2018) developed an optimization procedure to improve the designing and planning of multi waste biomass based supply chains according to CE for producing manifold kinds of bio-products through several kinds of technology in the same supply chain by the integration of mathematical modeling and fuzzy multi-objective decision-making.

**Supply chain and industry 4.0(I 4.0)**

, There are scare publications on smart manufacturing supply chain (SMSC) in comparison to plenty of investigations on SMSC according to converging and integrating ICPT such as Industry 4.0. A characteristic notion of convergence is Industry 4.0 by the German government. Industry 4.0—named the fourth industrial revolution—aiming at making industrial advancement toward a smart factory fortified with networked manufacture constituents and an effective manufacture system. With starting next-generation manufacture systems in Germany, a similar notion and strategies have been adopted by many countries around the globe to ensure their fabricating competitive ability. Recently, Industry 4.0 has been researched in domains concerning relevant government strategies, company policies, industrial usages, and the technology itself. In the majority of research on smart factories, the emphasis is of linking data and the cooperation of constituents and customizing products for market through ICPT.

A number of publications on a smart supply chain (SSC) are available in 'the smarter supply chain' of IBM [22], Cooke [23], and Dassault Systems [24]. Cooke [23] described the tendencies of the SSC including customer partitioning and omnichannel for sustainable business, and Dassault Systems [24] represented the SSC as achieving main gains of the supply chain including the quicker time-to-market and flexibility of re-optimizing plans. A number of investigators recommended enablers for providing enhanced ability to the supply chain. Tachizawa et al. [25] claimed that smart city and smart technologies furnished the substructure for the SSC. Dawid et al. [26] introduced manufacture procedures in which clients can partake. A commonly observed future supply chain is that all constituents have interconnections, decisions are taken independently, and flexibleness is necessary for responding to environmental variations.

The majority of planning models in supply chain consider minimalizing the overall expense of supply chain [27], and minimizing expense and environmental parameters taking sustainability into consideration [28–32]. Investigations available on optimization comprise gain maximization (or cost minimization) and time minimization—as in this research—. Melachrinoudis and Min [33] proposed maximization of the overall gain and aggregative location motivation for amenities and minimization of overall transportation time for determining optimum quantity and location. Two parameters, viz. minimization of the whole expenses and entire lateness of cycle time for reverse logistics, were examined by Du and Evans [34]. Liu and Papageorgiou [35] offered the production, distribution, and capacity planning of a universal supply chain for a planning model taking minimizing overall expenses, flow time, and lost sales into account. To discover an optimum planning approach, Tzeng et al. [36] and Liang [37] put forward multi-purposes of total cost and lead time. For a supply chain with flexibility [38–41], scarce cases of minimizing lead time are available. Moreover, the literature considered factors of profit or lead time to be attained in order. Nevertheless, a performance measure, SSCP, is proposed here the result of which combines two factors for careful management. Smart Supply Chain Performance (SSCP) puts forward a profit rate involving a tradeoff between profit and lead time, which have direct proportionality to quantity (profit =unit profit × quantity and lead time=unit required time × quantity). SSCP is an estimate of profit rate, then playing a role differently in SMSC under instantaneous data gathering. It is also capable of providing a lucrative policy with gain maximization and a reactive performace by minimization of lead time.

**Smart manufacturing supply chain**

The 'smart' means the capability of an object in recognizing per se the goal aimed to achieve, understanding and working well, and flexible responding to any variation. The integration and convergence of ICPT is necessary in the manufacture supply chain to have such a 'smart' supply chain as SMSC. Integration of ICPT provides SMSC supply chain with intelligence and flexibility. The objectives of SMSC correspond to the traditional supply chain: maintaining the balance of demand and supply, and maximizing customer gratification with minimal cost. Nonetheless, these objectives are achieved differently depending on application of ICPT in the supply chain. The entire constituents in the SMSC have connections on a communication network and collaborations on the custom-made products and services via information and communication technologies. Besides, they are able to be multi-functional through production technologies. Accordingly, the SMSC comprises three key characteristics of a connection, a collaboration, and a customization [42] and one attribute of flexible ability, which is a functionally and structurally relevant attribute in the SMSC.

The characteristics of SMSC are fundamental bases in forming the future supply chain. Connection and collaboration are supply-directed characteristics facilitating the response of the supply chain to variations.

The following describes every characteristic:

Among the above characteristics "Connection", emphasizes virtual connections between the whole members and constituents—wired or wireless—via information and communication networks. Such a connection allows for sharing and communicating by the supply chain in all information required from a real-time base, such that it can turn into a resource with intelligence, multi-functionality, informativity, and shareability per se [43, 44].

The other attribute, Collaboration, refers to the capability for communication, co-operation, and coordination with one another. As SMSC is regarded as a framework with virtual connection, more effective performance is obtained for physical and cyber constituents amongst the associates of the SMSC as a result of collaboration for solving problems emerging in the supply chain.

Additionally, "Customization" is a motivating driver of competitiveness as the major objective. Customization is applied to product, supply, dissemination channel, retailing, and buying methods. This feature makes producers capable of producing more custom-made products at inexpensive costs in a shorter lead time [47,48]

Altogether, the 'flexibility' of the supply chain is enhanced by the three features of SMSC [49], in particular, when congregated with more cutting-edge intelligent technologies. Flexibility denotes the capability in quick adaptation to the occurrence of difficulties [50] and is evidently the result of the three aforesaid characteristics, which can be classified into functional and structural flexibilities as described below.

### 3.3.1 Functional flexibility

**Volume and lead time flexibilities**: SMSC is constantly capable of controlling product type, amount, and lead time in an alternative fashion. Otherwise stated, according to Sabri and Beamon [38], Esmaeilikia et al. [50], and Rabbani et al. [51], the volume flexibility is the capability of changing the product quantity for a buffer. The buffer is the difference between the maximal and upper bound of capacity. Moreover, the lead time flexibility is where the SMSC has a delivery buffer, which is the capability of changing the delivered quantity and delivery date. The SMSC, with volume and time buffer as surplus, can react to unexpected fluctuations and fulfil production supplies other than the operation scheduled in advance.

**Capacity flexibility**: Capacity is the maximal permissible (production or distribution) amount per planning timeframe. Representative manageable capacities are the workforce, machine and facility, transport vehicle, among others. The tactical planning model (TPM) is capable of changing the maximal capacity through capacity flexibility [50–54].

**Decision-making flexibility**: Decision-making in SMSC occurs in real-time base onset. The introduced TPM supply is activated whenever existing supply chain experiences great fluctuations in main factors of the supply chain including production volume, capacity, lead time, etc. Because the whole data are gathered in real-time, SMSC undergoes continuous monitoring and managed with a decision maker who is able to implement the planning model or monitor the system performance with no human intervention.

### Structural flexibility

**Supplier flexibility:** It is one of decisive factors for network flexibility. Indeed, every supply chain component has its own supplier choices via selection of the most advantageous suppliers for materials or parts [50,51,54,56]. Since the SMSC possesses more custom-made and individualized product marketplaces, the producers require various extra parts and materials for even similar product type that needs more diverse supplier types.

**Channel flexibility**: Furthermore, Channel flexibility can assist clients in choosing their own preferred channels from the omnichannel including offline and online by channel flexibility. Using the omnichannel, customers are able to select their own delivery mode very cost-effectively and conveniently [57,58]. The customer can also receive service from different sources of producer, supplier, and seller in the form of raw materials, parts, modules, and finished goods [59,60].

According to the above features and traits, Figure 1 displays a basic schematic for the SMSC, in which the features reveal the philosophy. The supply chain tries for connection and collaboration with each other to become smart and offers customization as well.

The aforesaid flexibilities of the SMSC have associations with characteristics of the SMSC: connection, collaboration, and customization. Every kind of flexibility is revealed in the decision-making flow, TPM logic, and mathematical model.

**Problem definition**

Along with rising world population from 7 billion (2018) to 9 billion (2050), about 180 billion tones materials will be required yearly by 2050. Importantly, a plethora of these materials have been formulated, engineered and fabricated to be lifelong, disposable only, and not produced for reprocess or recycle. They function within the rules of a linear economy that takes, makes, uses and abuses materials, which at end of life (EoL) expire as wastage, i.e., crib to tomb. In societies, a huge bulk of wastage is generated as illustrated by waste plastics in recent years, which is obviously noticeable and in the purview of people. Currently, there is a need for changing this pattern and moving from linear to circular economies that consider resource, environment and society.

A pseudo real-time decision-making planning model is proposed in the present research. In essence, it is problematic to modify an existing supply plan when the operation and process are working in advance. Put differently, supply chain is not capable of adapting to the adjusted plan resulted from TPM throughout the lead time. As an example, a given truck has supplied 50 product units to the destination. It is not possible to employ a recently adapted amount in the course of the implementation time, despite the fact that the TPM leads to the supply of 100 products. Thus, supply chain can have time for changing and applying the re-planning solution in LTl. Hence, minimization of LTl creates short interval of decision-making with the aim of realizing pseudo real-time decision-making. The interval time shrinks between the first decision-making at the onset of the first $LT_{ij}$ (Fig. 3) and the second decision-making at the start of the second $LT_{ij}$, thereby securing the possibility of an approximate pseudo real-time decision-making. Figure 4 illustrates that SMSC can react to the varying environment by minimalizing LTl.

To obtain an efficient and stable SMSC at a definite level of performance, each parameter for the volume and lead time flexible abilities as buffer should be provided with a specified target or upper bound ($UU^P$ and $LT\ lr\ UP^{··}$). The upper bound in the SMSC contributes in preventing frequent planning, i.e. 'cover the changed data' in Figure 4, as frequent planning brings about extra costs and time. The upper bound of usage also has a contribution to the prevention of inventory and cycle time explosion due to a rise of usage [67]. The upper bound of lead time avoids the excess tardiness of product supply.

The actions of the SMSC in responding fluctuations in the environment consist of (1) the relocation of alternate routes by logistics flexibility [50,51,54,59], (2) the control of capacity by capacity flexibility [51–54], e.g. capacity development, overtime, and extra shifts, (3) altering the supply chain components by channel flexibility [57–60], and (4) the control of delivery quantity and lead time.

To summarize, the TPM model in Section 4 is tailored for the present and modified data with the aim of satisfying the director's decision. Once the present usage and lead time surpass the upper bound and the planning horizon ends, the supply chain director implements the TPM based on the real-time data gathered in the SMSC for effective response to the varying situation of the supply chain and obtaining novel supply chain plans. Otherwise, it is not possible to implement the TPM.